\documentclass[a4paper,reqno]{amsart}

\pdfoutput=1

\usepackage{nicefrac,amsmath,amssymb,amsthm,url,verbatim,enumerate,stmaryrd,mathtools,microtype}
\usepackage[utf8]{inputenc}
\usepackage[T1]{fontenc}
\usepackage{libertine}
\usepackage{dsfont}
\usepackage[libertine,cmintegrals,cmbraces]{newtxmath}
\usepackage[dvipsnames,svgnames,x11names]{xcolor}
\definecolor{darkred}{RGB}{139,0,0}
\definecolor{darkblue}{RGB}{0,0,139}
\definecolor{darkgreen}{RGB}{0,100,0}
\usepackage[mathscr]{euscript}
\usepackage[pagebackref]{hyperref}
\hypersetup{
	colorlinks   = true, 
	urlcolor     = darkred, 
	linkcolor    = darkred, 
	citecolor   = darkgreen}
\usepackage{tikz}
\usepackage{tikz-cd}
\usepackage{pgfplots}
\usepackage{calc}
\usepgfmodule{nonlineartransformations}
\usepgflibrary{curvilinear}
\usetikzlibrary{decorations.markings,patterns,arrows.meta}
\tikzset{>={Straight Barb[scale=0.8]}}
\tikzcdset{
	arrow style=tikz,
	diagrams={>={Straight Barb[scale=0.8]}}
}
\usepackage[capitalise,noabbrev]{cleveref}

\calclayout

\AtBeginDocument{%
   \def\MR#1{}
}

\newcommand\smallsquare{{\mathbin{\text{\raise0.17ex\hbox{\scalebox{.7}{$\blacksquare$}}}}}}

\newcommand{\BO}{\ensuremath{\mathrm{BO}}}
\newcommand{\BSO}{\ensuremath{\mathrm{BSO}}}
\newcommand{\Fr}{\ensuremath{\mathrm{Fr}}}
\newcommand{\Diff}{{\ensuremath{\mathrm{Diff}}}}

\newcommand{\BDiff}{\ensuremath{\mathrm{BDiff}}}

\newcommand{\Homeo}{\ensuremath{\mathrm{Homeo}}}

\newcommand{\BHomeo}{\ensuremath{\mathrm{BHomeo}}}

\newcommand{\Map}{\ensuremath{\mathrm{Map}}}

\newcommand{\Emb}{\ensuremath{\mathrm{Emb}}}

\newcommand{\BEmb}{\ensuremath{\mathrm{BEmb}}}
\newcommand{\BAut}{\ensuremath{\mathrm{BAut}}}
\newcommand{\BSAut}{\ensuremath{\mathrm{BSAut}}}

\newcommand{\BTop}{\ensuremath{\mathrm{BTop}}}
\newcommand{\BSTop}{\ensuremath{\mathrm{BSTop}}}

\newcommand{\inc}{\ensuremath{\mathrm{inc}}}

\newcommand{\Fun}{\ensuremath{\mathrm{Fun}}}

\newcommand{\PSh}{\mathrm{PSh}}

\DeclareMathAlphabet{\mathpzc}{OT1}{pzc}{m}{it}

\renewcommand{\rm}[1]{\ensuremath{\mathrm{#1}}}
\newcommand{\icat}[1]{\ensuremath{\mathscr{#1}}}

\newcommand{\bfR}{\ensuremath{\mathbf{R}}}
\newcommand{\bfZ}{\ensuremath{\mathbf{Z}}}
\newcommand{\bfQ}{\ensuremath{\mathbf{Q}}}
\newcommand{\bfS}{\ensuremath{\mathbf{S}}}

\newcommand{\Cosp}{\ensuremath{\mathrm{Cosp}}}

\newcommand{\colim}{\ensuremath{\mathrm{colim}}}

\newcommand{\lra}{\longrightarrow}

\newcommand{\scr}[1]{\mathscr{#1}}
\newcommand{\ul}[1]{\underline{#1}}

\newcommand{\GL}{\mathrm{GL}}

\newcommand{\Aut}{\mathrm{Aut}}

\newcommand{\coker}{\mathrm{coker}}
\newcommand{\bP}{\mathrm{bP}}

\renewcommand{\Top}{\mathrm{Top}}

\newcommand{\id}{\mathrm{id}}

\newcommand{\fib}{\mathrm{fib}}

\newcommand{\op}{\mathrm{op}}

\newcommand{\Pro}{\mathrm{Pro}}

\newcommand{\un}{\mathrm{un}}
\newcommand{\red}{\mathrm{red}}

\newcommand{\col}{\mathrm{col}}

\newcommand{\Env}{\ensuremath{\mathrm{Env}}}

\newcommand{\Th}{{\mathrm{Th}}}

\newcommand{\fQ}{{\mathrm{f}\bfQ}}

\newcommand{\Disc}{\ensuremath{\rm{Disc}}}
\newcommand{\DiscInf}{\ensuremath{\icat{D}\mathrm{isc}}}

\newcommand{\ManInf}{\ensuremath{\icat{M}\mathrm{an}}}
\newcommand{\Fin}{\ensuremath{\mathrm{Fin}}}
\newcommand{\ALG}{\ensuremath{\mathrm{ALG}}}

\newcommand{\ncBordInf}{\ensuremath{\mathrm{nc}\icat{B}\mathrm{ord}}}

\newcommand{\Mul}{\ensuremath{\mathrm{Mul}}}

\newcommand{\oo}{\ensuremath{\mathrm{o}}}
\newcommand{\ot}{\ensuremath{\mathrm{t}}}
\newcommand{\Spc}{\ensuremath{\mathscr{S}\rm{pc}}}

\newtheorem{thm}{Theorem}[section]

\newtheorem{question}[thm]{Question}
\newtheorem{prob}[thm]{Problem}

\theoremstyle{definition}
\newtheorem{dfn}[thm]{Definition}

\newtheorem*{nconvention}{Convention}
\theoremstyle{remark}
\newtheorem{ex}[thm]{Example}

\newtheorem{rem}[thm]{Remark}

\begin{document}

\title{Manifolds and Disc-presheaves}

\author{Alexander Kupers}
\address{Department of Computer and Mathematical Sciences, University of Toronto Scarborough, 1265 Military Trail, Toronto, ON M1C 1A4, Canada}
\email{a.kupers@utoronto.ca}

\begin{abstract}This essay explains an approach to the study of smooth manifolds which compares them to presheaves on a category of discs, also known as embedding calculus. We highlight recent work that shows this approach has many desirable properties, as well as recent applications demonstrating its strength.\end{abstract}

\vspace*{-1cm}
\maketitle

\vspace{-.5cm}
\tableofcontents

\vspace{-.5cm} A guiding problem of algebraic topology is the study of manifolds and their automorphisms, asking for a description of moduli spaces of manifolds. One strategy for this is to compare the category of manifolds to one satisfying the following desiderata:
\begin{enumerate}[(I)]
	\item \label{enum:des-prop} It has analogous properties and structure as the category of manifolds.
	\item \label{enum:des-cat} It can itself be understood.
	\item \label{enum:des-sdisc} Its difference to the category of manifolds can be understood.
\end{enumerate}
More concretely, Desideratum \eqref{enum:des-prop} asks that we can mimic the constructions we perform on manifolds, so we can study these constructions as well as use them in our proofs. Desiderata \eqref{enum:des-cat} and \eqref{enum:des-sdisc} ask for descriptions in terms of computationally accessible homotopy theory and algebra. The success of such a strategy can be measured by its applications and how it spurs progress in related fields.

\smallskip

This essay is an extended argument in favour of comparing manifolds to \emph{Disc-presheaves}. We first explain how we have learned that it satisfies the desiderata \eqref{enum:des-prop}--\eqref{enum:des-sdisc}, at least to a large extent \cite{KKSDisc,KKOperadic}. We then illustrate how it is successful by the above standards. Firstly, Disc-presheaves are a main ingredient---in addition to the paired techniques of homological stability and stable moduli spaces---in many of the recent breakthroughs in the subject, among others, on the rational homotopy of diffeomorphisms of discs \cite{KRWDiscs,KKOdd,RandalWilliamsDiscs} and on Pontryagin--Weiss classes \cite{WeissPontryagin,KKSquare}. Secondly, the category of Disc-presheaves is closely connected to operads and graph complexes, and recent work suggests links to algebraic $K$- and $L$-theory. Finally, it has several advantages over the classical approach to manifold theory using surgery theory and pseudoisotopy theory.


\begin{nconvention}This survey is written in the language of $\infty$-categories. A reader unfamiliar with this will not lose much by substituting categories enriched in spaces for $\infty$-categories.
\end{nconvention}

\subsection*{Acknowledgments} As will be clear to the reader, the work described in this essay is deeply inspired to the groundbreaking work of Michael Weiss. It is also a pleasure to acknowledge my collaborators on the subject of this essay: Mauricio Bustamante, Manuel Krannich, Ben Knudsen, Fadi Mezher, and Oscar Randal-Williams. I also thank S\o ren Galatius, Manuel Krannich, Oscar Randal-Williams, and Ismael Sierra for feedback on an earlier version. AK acknowledges the support of the Natural Sciences and Engineering Research Council of Canada (NSERC) [funding reference number 512156 and 512250]. 

\section{Disc-presheaves}\label{sec:motivation} In this first section we introduce and motivate the subject of this essay---the Disc-presheaves of \cref{def:disc-presheaf}. This concept first arose in embedding calculus, which can be interpreted as studying mapping spaces of Disc-presheaves associated to manifolds, see \cref{sec:discpresheaves}.

\subsection{A pointillist view of manifolds} How do we distinguish smooth manifolds from each other, up to diffeomorphism? The usual invariants from algebraic topology are homotopy-invariant, and thus only detect the underlying homotopy type. However, we can expand their scope by constructing other topological spaces from a smooth manifold $M$ (possibly with some additional data) and distinguishing their homotopy types instead. One could take an analytic approach and consider for instance spaces of solutions to systems of partial differential equations, leading to gauge-theoretic invariants. We rather take a geometric approach and consider spaces of smooth embeddings of very simple manifolds into $M$.

\subsubsection{Embeddings of points} \label{sec:emb-pts} The simplest manifolds are arguably collections of points: $\ul{k} \coloneqq \{1,\ldots,k\}$. We let $\Emb^\oo(-,-)$ denote spaces of smooth embeddings in the smooth topology.

\begin{dfn}The \emph{ordered configuration space of $k$ points} is given by $\Emb^\oo(\ul{k},M)$.
\end{dfn}

The homotopy types of ordered configuration spaces are diffeomorphism-invariants of $M$, but it is a priori not clear they are \emph{not} homotopy-invariants. However, there is the following example:

\begin{ex}The 3-dimensional lens spaces $L(7,1)$ and $L(7,2)$ are not diffeomorphic but \emph{are} homotopy equivalent, and Longoni and Salvatore showed that $\smash{\Emb^\oo}(\ul{2},L(7,1))$ is not homotopy equivalent to $\smash{\Emb^\oo}(\ul{2},L(7,2))$ \cite{LongoniSalvatore}.\end{ex}

A defect of the ordered configuration spaces is that their homotopy types \emph{are} homeomorphism-invariants, as $\Emb^\oo(\ul{k},M)$ is homeomorphic to $\{(m_1,\ldots,m_k) \in M^k \mid m_i \neq m_j \text{ if } i \neq j\}$.

\subsubsection{Embedding of open discs} \label{sec:emb-open-discs} To recover the smooth structure, at least we need to record the tangent bundle. If $M$ is $d$-dimensional, this can be done by replacing the finite collections of points $\ul{k}$ with the finite collections of open discs $\ul{k} \times \bfR^d$. In the case $k=1$, taking the total derivative at the origin of $\bfR^d$ yields an equivalence $\smash{\Emb^\oo(\bfR^d,M) \xrightarrow{\simeq} \Fr(TM)}$ to the total space of frame bundle of the tangent bundle $TM$. For general $k$, embedding spaces of open discs are equivalent to \emph{framed ordered configuration spaces}:
\[\Emb^\oo(\ul{k} \times \bfR^d,M) \overset{\simeq}\lra \Emb^\oo(\ul{k},M) \times_{M^{\ul{k}}} \Fr(TM)^{\ul{k}}.\]
To obtain not just the total space $\Fr(TM)$ but its structure as a principal $\GL_d(\bfR)$-bundle---equivalently, the structure of $TM$ as a $d$-dimensional vector bundle---we can make use of naturality in the domain: the inclusion $\GL_d(\bfR) \to \Emb^\oo(\bfR^d,\bfR^d)$ of linear maps into self-embeddings gives an $\GL_d(\bfR)$-action on $\Emb^\oo(\bfR^d,M)$ that corresponds to the action on $\Fr(TM)$ by reframing.

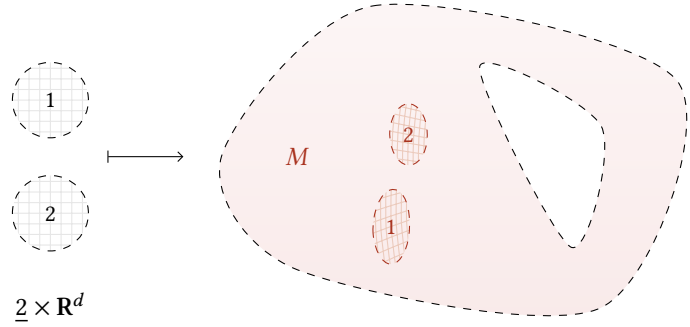
\begin{figure}
\begin{tikzpicture}
	\draw [dashed,shade,top color=Mahogany!3!white,
	bottom color=Mahogany!7!white] plot [smooth cycle,tension=.5,scale=2] coordinates {(0, 0) (1,1) (3,0.5) (2.5,-1) (.5,-.75)};
	\draw [dashed,fill=white] plot [smooth cycle,tension=.5,scale=2] coordinates {(1.7,.6) (2.5,.2) (2.3,-.6)};
	\node at (1,0) [Mahogany] {$M$};
	
	\begin{scope}[shift={(-2.25,.75)}]
	\begin{scope}
		\clip (0,0) circle (.5cm);
		\draw [thin,step=4pt,color=black!10!white] (-2,-2) grid (2,2);
	\end{scope}
	\draw[dashed,black] (0,0) circle (.5cm);
	\node at (0,0) {\small $1$};
	\end{scope}
	
	\begin{scope}[shift={(-2.25,-.75)}]
		\begin{scope}
			\clip (0,0) circle (.5cm);
			\draw [thin,step=4pt,color=black!10!white] (-2,-2) grid (2,2);
		\end{scope}
		\draw[dashed,black] (0,0) circle (.5cm);
		\node at (0,0) {\small $2$};
	\end{scope}
	
	\node at (-2.25,-2) {$\ul{2} \times \bfR^d$};
	
	\draw [|->] (-1.5,0) -- (-.5,0);
	
	\begin{scope}[shift={(2.5,.8)},rotate=90]
	\scoped{
		\pgfsetcurvilinearbeziercurve
		{\pgfpointxy{5}{-1}}{\pgfpointxy{5}{9}}
		{\pgfpointxy{40}{15}}{\pgfpointxy{25}{20}}
		\pgftransformnonlinear{\pgfgetlastxy\x\y%
			\pgfpointcurvilinearbezierorthogonal{\y}{\x+.2*\y}}%
		\begin{scope}
			\clip (0,0) circle (.25cm);
			\draw [thin,step=2pt,color=Mahogany!20!white] (-2,-2) grid (2,2);
				\node at (0,0) [Mahogany] {\small 2};
		\end{scope}
		\draw[dashed,Mahogany] (0,0) circle (.25cm);
	}
	\end{scope}
	
	\begin{scope}[shift={(2.8,-.2)}]
		\scoped{
			\pgfsetcurvilinearbeziercurve
			{\pgfpointxy{5}{0}}{\pgfpointxy{5}{10}}
			{\pgfpointxy{55}{15}}{\pgfpointxy{15}{20}}
			\pgftransformnonlinear{\pgfgetlastxy\x\y%
				\pgfpointcurvilinearbezierorthogonal{\y-.1*\x}{\x-.2*\y}}%
			\begin{scope}
				\clip (0,0) circle (.25cm);
				\draw [thin,step=2pt,color=Mahogany!20!white] (-2,-2) grid (2,2);
				\node at (0,0) [Mahogany] {\small 1};
			\end{scope}
			\draw[dashed,Mahogany] (0,0) circle (.25cm);
		}
	\end{scope}
\end{tikzpicture}
\caption{An element of $\Emb^\oo(\ul{2} \times \bfR^d,M)$. The map to the framed ordered configuration space records the centres of the two open discs on the right, as well as the induced frame in the tangent spaces at those points.}
\end{figure}

\subsubsection{Embeddings of open discs and maps between them} \label{sec:emb-maps} There are further maps that can be constructed from the naturality in the domain. Indeed, \emph{any} embedding $\ul{k}' \times \bfR^d \hookrightarrow \ul{k} \times \bfR^d$ gives rise to a map
\[\Emb^\oo(\ul{k} \times \bfR^d,M) \lra \Emb^\oo(\ul{k}' \times \bfR^d,M),\]
which depends continuously on the embedding. Taking it to be a permutation of components yields a map that permutes points, taking it to be an inclusion of components a map that forgets some points, and taking it to be an inclusion of several open discs into a single one a map that splits a point.

\subsection{Disc-presheaves} Boavida de Brito and Weiss observed that the above data---the spaces $\Emb^\oo(\ul{k} \times \bfR^d,M)$ for $k \geq 0$ and maps between them from naturality in the domain---can be encoded in a single mathematical object, which is the main subject of this survey \cite{BoavidadeBritoWeissSheaves}.

Let $\ManInf^\oo_d$ be the category whose objects are (possibly non-compact) $d$-dimensional smooth manifolds and whose mapping spaces are spaces of smooth embeddings. This has a full subcategory $\DiscInf^\oo_d$ on those objects diffeomorphic to $S \times \bfR^d$ for $S$ a finite set. For any category $\scr{C}$ we can form a category $\PSh(\scr{C}) \coloneqq \Fun(\scr{C}^\op,\Spc)$ of presheaves valued in the category $\Spc$ of spaces.

\begin{dfn}\label{def:disc-presheaf} We refer to an object $X \in \PSh(\DiscInf_d^\oo)$ as a \emph{Disc-presheaf}.\end{dfn}

Unwinding the definition, a Disc-presheaf $X$ is a collection of spaces $X(\ul{k} \times \bfR^d)$ for $k \geq 0$ that are contravariantly natural in embeddings of the input. The Disc-presheaf given by
\[\begin{aligned} E_M \colon \DiscInf_d^\op &\lra \Spc \\
	\ul{k} \times \bfR^d &\longmapsto \Emb^\oo(\ul{k} \times \bfR^d,M)\end{aligned}\] 
is the image of $M$ under the restricted Yoneda embedding
\begin{equation}\label{eqn:efunctor}E \colon \ManInf^\oo_d \overset{y}\lra \PSh(\ManInf^\oo_d) \overset{\iota^*}\lra \PSh(\DiscInf^\oo_d)\end{equation}
along the inclusion $\iota \colon \DiscInf^\oo_d \hookrightarrow \ManInf_d^\oo$.

\begin{dfn}We refer to $E_M \in \PSh(\DiscInf_d^\oo)$ as the \emph{Disc-presheaf associated to the $d$-dimensional smooth manifold $M$}.\end{dfn}

Since $E$ is a functor, it induces maps of spaces $\Emb^\oo(M,M') \to \smash{\Map_{\PSh(\DiscInf_d^\oo)}(E_M,E_{M'})}$ that are compatible with composition. Thus diffeomorphisms are sent to invertible maps and we get a map $\Diff(M) \to \Aut_{\PSh(\DiscInf_d^\oo)}(E_M)$. Consequently, not only can Disc-presheaves be used to study smooth manifolds, but the induced maps on Disc-presheaves give information about spaces of smooth embeddings or diffeomorphisms.

\subsection{Moduli spaces and structure spaces} \label{sec:moduli-sdisc} So far our stated motivation was to find invariants that distinguish manifolds, but it is equally interesting to find manifolds that realise those invariants: this means constructing manifolds or diffeomorphisms, realising Disc-presheaves or automorphisms of these. To capture these goals, we use the language of moduli spaces and structure spaces. The full subcategory $\ManInf_d^\rm{cl}$ of $\ManInf_d$ consisting of closed $d$-dimensional smooth manifolds is an $\infty$-groupoid---or equivalently, a space---as embeddings between such manifolds must be diffeomorphisms, and $E$ maps this to the $\infty$-groupoid core $\PSh(\DiscInf^\oo_d)^\simeq$. This gives a map
\begin{equation}\label{eqn:emap} E \colon \ManInf_d^\rm{cl} \lra \PSh(\DiscInf^\oo_d)^\simeq\end{equation}
from the moduli space of closed $d$-dimensional smooth manifolds to that of Disc-presheaves. From this perspective, Desideratum \eqref{enum:des-cat} seeks an answer to the following question, or more generally for a description of the category $\PSh(\DiscInf^\oo_d)$:

\begin{question}\label{ques:disc-presheaves} What is the homotopy type of the moduli spaces $\PSh(\DiscInf^\oo_d)^\simeq$?
\end{question}

Given such a description, one can use Disc-presheaves to study manifolds by answering the following question:

\begin{question}\label{ques:moduli} How close is the map \eqref{eqn:emap} to being an equivalence?\end{question}

We can phrase an answer in terms of the difference between the domain and target of map \eqref{eqn:emap}. This is captured by the following definition from \cite{KKSDisc}, in analogy with surgery theory:

\begin{dfn}The \emph{Disc-structure space} of $X \in \PSh(\DiscInf_d^\oo)$ is the fibre
	\[S^\Disc(X) \coloneqq \fib_X\big[\ManInf_d^\rm{cl} \overset{E}\lra \PSh(\DiscInf^\oo_d)^\simeq\big].\]
\end{dfn}

The Disc-structure space $S^\Disc(X)$ can be thought of as the space of ``manifold structures on $X$''. The set of path-components is in in bijection with equivalences $X \simeq E_M$ up to maps induced by diffeomorphisms on the right side, and such a path component is equivalent to the component of $\Aut_{\PSh(\DiscInf^\oo_d)}(X)/\Diff(M)$ containing the identity. From this perspective, Desideratum \eqref{enum:des-sdisc} asks for an answer to the following question:

\begin{question}\label{ques:sdisc} What is the homotopy type of the Disc-structure space $S^{Disc}(X)$?\end{question}

There has been significant recent progress on these questions, suggesting that Disc-presheaves sit at a convenient midpoint between smooth manifolds and homotopy types. More precisely, the answer to \cref{ques:moduli} is generally \emph{no} (\cref{thm:non-triviality} and \ref{thm:2-type-invariance}) but what we know about \cref{ques:sdisc} indicates that the failure of \eqref{eqn:emap} is captured by a highly-structured mathematical object (\cref{thm:infinite-loop-space}). It suggests a future in which one can study manifolds through the homotopy theory of Disc-presheaves and Disc-structure spaces. 

\section{Desideratum \eqref{enum:des-prop}: properties and structure of Disc-presheaves} \label{sec:structure} To make use of Disc-presheaves in manifold theory, we need analogues of such constructions as disjoint union, gluing, cartesian products, covering spaces, and tangential structures, as well as of such results as isotopy extension and smoothing theory. These are provided by the foundations of \cite{KKSDisc,KKOperadic} developed by Krannich and myself, which connects neatly to the theory of operads; these papers serve as the reference for this section unless mentioned otherwise.

\subsection{A field-theoretic structure} We start by explaining how the category of Disc-presheaves admits disjoint union and gluing operations, encoded in a field-theoretic analogous to Atiyah's approach via symmetric monoidal functors from bordism categories.

\subsubsection{Disjoint union and gluing by algebras and modules} \label{sec:disjoint-union-gluing} 
An embedding $S \times \bfR^d \to M_1 \sqcup M_2$ into a disjoint union is uniquely determined by a map $\psi \colon S \to \smash{\ul{2}}$ of finite sets and a pair of embeddings $\psi^{-1}(i) \times \bfR^d \to M_i$ for $i=1,2$. This exhibits the Disc-presheaf of a disjoint union of two manifolds as the Day convolution with respect to disjoint union of discs
\[E_{M_1 \sqcup M_2} \simeq E_{M_1} \otimes E_{M_2} \text{ where } (E_{M_1} \otimes E_{M_2})(S \times \bfR^d) \simeq \bigsqcup_{S_1\sqcup S_2 = S} E_{M_1}(S_1 \times \bfR^d) \times E_{M_2}(S_2 \times \bfR^d),\]
lifting $E$ to a symmetric monoidal functor. For $(d-1)$-dimensional smooth manifold $Q$, the product $Q \times \bfR$ admits the structure of $E_1$-algebra in $\ManInf_d$ with multiplication given by
\[Q \times \bfR \sqcup Q \times \bfR\lra Q \times \bfR\]
induced by an embedding $\bfR \sqcup \bfR \to \bfR$. Similarly, a choice of collar for $Q \subset \partial M$ makes $M$ into a left $E_1$-module over $Q \times \bfR$. Applying $E$, this induces on $E_M$ the structure of a left $E_1$-module over $E_{Q \times \bfR}$. The Disc-presheaf of a gluing of two manifolds along a common codimension zero submanifold of their boundary is given in these terms by the relative tensor product
 \[E_{M_1 \sqcup_Q M_2} \simeq E_{M_1} \otimes_{E_{Q \times \bfR}} E_{M_2}.\]

The Disc-presheaves arising from manifolds have the property that their value on $\varnothing$ is contractible; we say they are \emph{unital}. These are preserved by (relative) tensor products, and the coherence properties of these constructions are encoded in the existence of a symmetric monoidal $(\infty,2)$-category $\ALG(\PSh^\rm{un}(\DiscInf_d^\oo))$. An $(\infty,2)$-category is informally a category enriched in categories: it thus has a category of objects (here given by the category of $E_1$-algebras in the category $\PSh^\rm{un}(\DiscInf_d^\oo)$ of unital Disc-presheaves) as well as mapping categories between objects (here the mapping category from $A$ to $B$ is the category of $A,B$-bimodules in unital Disc-presheaves) with composition functors (here given by relative tensor products). The additional symmetric monoidal structure is given by external tensor products.

\subsubsection{From bordisms to bimodules}
By design, the ``disjoint union'' and ``gluing'' constructions for Disc-presheaves are such that the construction of the restricted Yoneda functor $E \colon \ManInf_d \to \PSh^\rm{un}(\DiscInf_d^\oo)$ lifts to a functor of symmetric monoidal $(\infty,2)$-categories
\begin{equation}\label{eqn:e-2cat} E \colon \rm{ncBord}_d^\oo \lra \ALG(\PSh^\rm{un}(\DiscInf_d^\oo)).\end{equation}
The domain is a non-compact bordism category. Its category of objects has as objects possibly non-compact $(d-1)$-dimensional manifolds and as mapping spaces embedding spaces. It has a mapping category from $P$ to $Q$, which has as objects possibly non-compact bordisms $W \colon P \leadsto Q$ and as mapping spaces embedding spaces relative to $P$ and $Q$. Composition is induced by gluing and it is symmetric monoidal under disjoint union. If $M$ and $N$ have common boundary $\partial M = P = \partial N$, taking $Q = \varnothing$ we get a map
\[\Emb_\partial(M,N) \lra \Map_{\rm{LMod}_{E_{P \times I}}(\PSh^\rm{un}(\DiscInf_d^\oo)}(E_M,E_N),\]
leading us to interpret the mapping space of left $E_1$-modules over $E_{P \times I}$ in $\PSh^\rm{un}(\DiscInf_d^\oo)$, as maps of presheaves relative to the boundary.

\begin{rem}One could interpret \eqref{eqn:e-2cat} as a field theory---closely related to \emph{factorisation homology}, see e.g.~\cite{AyalaFrancis}---though its domain differs from more familiar bordism categories by allowing non-compact rather than compact manifolds and embeddings rather than diffeomorphisms. \end{rem}

\subsection{Tangential structures and smoothing theory}\label{sec:smoothing-theory} It is often useful to endow manifolds with the following additional structure: a tangential structure is a map $\theta \colon B \to \BO(d)$ and a $\theta$-structure on a smooth $d$-dimensional manifold $M$ is a lift of a map $M \to \BO(d)$ classifying its tangent bundle along $\theta$. For example, for $\BSO(d) \to \BO(d)$ this is an orientation and for $\ast \to \BO(d)$ this is a framing. These structures appear in the comparison of smooth to topological manifolds: smoothing theory says that if $d \neq 4$ there is a pullback square
\begin{equation}\label{eqn:smoothing-theory} \begin{tikzcd} \ManInf_d^\oo \rar \dar & \Spc_{/\BO(d)} \dar \\[-7pt]
	\ManInf_d^\ot \rar & \Spc_{/\BTop(d)} \end{tikzcd}\end{equation}
where $	\ManInf_d^\ot$ is the category whose objects are $d$-dimensional topological manifolds and whose mapping spaces are spaces of topological embeddings \cite{KirbySiebenmann}. The left vertical map forgets the smooth structure, the right vertical map composes with $\BO(d) \to \BTop(d)$, and the horizontal maps extract the classifying map for the smooth or topological tangent bundle.

Replacing $\DiscInf_d^\oo \subset \ManInf_d^\oo$ with $\DiscInf_d^\ot \subset \ManInf_d^\ot$, we obtain a topological variant $E^\ot$ of \eqref{eqn:efunctor} with analogous properties. In \cref{sec:emb-open-discs} we saw that the map classifying the smooth tangent bundle can be recovered from $E_M$ as $M \simeq E_M(\bfR^d)_{\rm{O}(d)} \to \BO(d)$, where $(-)_{\rm{O}(d)}$ denotes (homotopy) orbits. It is similarly true, with $\Top(d) \coloneqq \Homeo(\bfR^d)$ taking the role of $\rm{O}(d)$, that the map classifying the topological tangent bundle can be recovered as $M \simeq \smash{E^\ot_M(\bfR^d)_{\rm{Top}(d)}} \to \BTop(d)$. Thus we can factor \eqref{eqn:smoothing-theory} over Disc-presheaves, compatibly with disjoint union and gluing, in the form of a commutative diagram
\[\begin{tikzcd}
	\ncBordInf^\oo_d \rar \dar & \ALG(\PSh^\rm{un}(\DiscInf_d^\oo)) \rar \dar & \Cosp(\Spc_{/\BO(d)}) \dar \\[-7pt]
	\ncBordInf^\ot_d \rar & \ALG(\PSh^\rm{un}(\DiscInf_d^\ot)) \rar & \Cosp(\Spc_{/\BTop(d)})
\end{tikzcd}\]
of symmetric monoidal $(\infty,2)$-categories, whose right terms are cospan categories. The right square is a pullback (``smoothing theory for Disc-presheaves'') and the outer square is a pullback on mapping spaces in mapping categories if $d \neq 4$, hence so is the left square.

\begin{ex}\label{rem:dim-4} The right square is a pullback for \emph{all} dimensions, so Disc-presheaves behave as if smoothing theory also applies in dimension $d=4$. That is, they are only sensitive to the \emph{formally smooth manifold} given by the underlying topological manifold and the lift of its topological tangent bundle to a vector bundle. This idea underlies \cite{KKnudsen}.
\end{ex}

\subsection{Operadic foundations} \label{sec:operadic} It is possible and convenient to reframe the previous discussion in terms of operads. 

\subsubsection{Disc-presheaves as right $E_d^\oo$-modules} We start with the observation, going back at least to work of Arone, Lambrechts, Turchin, and Volić \cite{ALTV}, that Disc-presheaves are right modules over the framed little $d$-discs operad. The framed little $d$-discs operad $\smash{E_d^\oo}$ has a single colour that we call $\bfR^d$, space of $n$-ary operations given by a space of smooth embeddings 
\[\Mul^{E_d^\oo}((\bfR^d)_{s \in S},\bfR^d) \simeq \Emb^\oo(\sqcup_{s \in S} \bfR^d,\bfR^d),\]
and operadic composition given by composition of embeddings. It has a symmetric monoidal envelope $\Env(E_d^\oo)$ whose objects are finite collection of colours $(\bfR^d)_{s \in S}$, morphism spaces built from multi-operations as 
\[\vspace{-.1cm} \Map_{\Env(E_d^\oo)}((\bfR^d)_{s \in S},(\bfR^d)_{t \in T}) = \bigsqcup_{\varphi \colon S \to T} \bigsqcap_{t \in T} \Mul^{E_d^\oo}((\bfR^d)_{s \in \varphi^{-1}(t)},\bfR^d),\]
and symmetric monoidal structure given by disjoint union. Since an embedding $\sqcup_{s \in S} \bfR^d \hookrightarrow \sqcup_{t \in T} \bfR^d$ is uniquely determined by a map $\varphi \colon S \to T$ of of finite sets and embeddings $\sqcup_{s \in \varphi^{-1}(t)} \bfR^d \to \bfR^d$ for all $t \in T$, we recognise there is an equivalence of symmetric monoidal categories $\DiscInf_d^\oo \simeq \Env(E_d^\oo)$. This induces an equivalence $\PSh(\DiscInf_d^\oo) \simeq \rm{rMod}(E_d^\oo)$, with right modules over an operad $\scr{O}$ defined as presheaves on its symmetric monoidal envelope: $\rm{rMod}(\scr{O}) \coloneqq \PSh(\Env(\scr{O}))$.

\subsubsection{The ``particle'' tangential structure} \label{sec:particle} The target of \eqref{eqn:e-2cat} makes sense for any operad $\scr{O}$---defining right $\scr{O}$-modules as presheaves on the symmetric monoidal envelope $\Env(\scr{O})$. In this context there exists a useful generalisation of smoothing theory for Disc-presheaves.

If all spaces of 0-ary multi-operations of an operad $\scr{O}$ are contractible we say it is \emph{unital} and if moreover 1-ary multi-operations are contractible we say it is \emph{reduced}. The colours and $1$-ary multi-operations assemble to a category of colours $\scr{O}^\col$, and by work of Lurie a unital operad $\scr{O}$ so that $\scr{O}^\col$ is a groupoid can be written as a colimit $\scr{O} \simeq \colim_{\scr{O}^\col} \scr{O}^\red$ for a reduced operad $\scr{O}^\red$ and a map $\theta \colon \scr{O}^\col \to \BAut(\scr{O}^\red)$. For example, $\smash{E_d^\oo \simeq \colim_{\BO(d)} E_d}$ for the action of $\rm{O}(d)$ on the little $d$-discs operad $E_d$ by rotation, and similarly for $E_d^\ot$. We think of a map $\theta \colon B \to \BAut(E_d)$ as a tangential structure for Disc-presheaves and for any map over $\BAut(E_d)$ there is a pullback square of symmetric monoidal $(\infty,2)$-categories
\[\vspace{-.1cm} \begin{tikzcd}
	\ALG(\rm{rMod}(E_d^\theta)) \rar \dar & \Cosp(\Spc_{/B}) \dar \\[-7pt]
	\ALG(\rm{rMod}(E_d^{\theta'})) \rar & \Cosp(\Spc_{/B'}).
\end{tikzcd}\]
This implies that the categories of Disc-presheaves obtained from the operads $\smash{E_d^\oo}$ or $\smash{E_d^\ot}$ are obtained by merely adding a tangential structure to the more fundamental category of Disc-presheaves obtained from the operad $E_d^\rm{p} \coloneqq \smash{\colim_{\BAut(E_d)} E_d}$; the superscript $p$ stands for ``particle''. We will this terminology by relating it in forthcoming work with Krannich to the configuration categories studied Boavida de Brito and Weiss \cite{BoavidadeBritoWeissConfCat}.

\begin{ex}\label{ex:alexander-trick} There is a delooping result relating Disc-presheaves with tangential structure and maps of operads: using that $E_{D^d}$ is a left $E_{S^{d-1} \times \bfR}$-module and adding tangential structures, it says
\vspace{-.2cm}\begin{equation}\label{eqn:delooping}\Aut_{\rm{LMod}_{E_{S^{d-1} \times \bfR}}(\rm{rMod}(E_d^\theta))}(E_{D^d}) \simeq \Omega^{d+1} \fib\big(\theta \colon B \to \BAut(E_d)\big),\end{equation}
a variant of a result of Dwyer and Hess \cite{DwyerHess}. For $\theta = \id_{\BAut(E_d)}$ this is contractible, an instance of the Alexander trick, c.f.~\cite{BoavidadeBritoWeissConfCat}.\end{ex}

\subsection{Outlook: further structure} We have not yet discussed isotopy extension, but a partial analogue exists in case of convergence for the base, see \cref{sec:isotopy-extension}. There is further structure on the category of manifolds whose analogues have not yet been developed for Disc-presheaves. We intend to in the future, when applications require them:
\begin{enumerate}[\noindent (1)]
	\item Extending \eqref{eqn:e-2cat} to fully extended field-theoretic structure that is a functor of symmetric monoidal $(\infty,d)$-categories.
	\item Relating cartesian products of manifolds to Boardman--Vogt tensor products of operads and right modules.
	\item Constructing variants for manifolds with boundary (or corners), and embeddings that preserve these setwise in terms of the swiss cheese operads (or their generalisations).
\end{enumerate}

\section{Desideratum \eqref{enum:des-cat}: understanding Disc-presheaves} \label{sec:discpresheaves}Given a well-behaved comparison between manifolds and Disc-presheaves, we now discuss how one may understand the latter.

\subsection{The cardinality tower} A powerful tool to study Disc-presheaves is a tower  obtained by filtering the category $\DiscInf_d^\oo$ by cardinality. This was the original motivation for studying them, in the setting of embedding calculus, as we recall now following \cite{KKOperadic}.

\subsubsection{The tower and its layers}  The \emph{cardinality filtration} of the category $\DiscInf_d^\oo$ by full subcategories is obtained bounding the number of open discs we allow:
\[\ast \simeq \DiscInf_{d,\leq 0}^{\oo} \overset{\subset}{\lra} \DiscInf_{d,\leq 1} \overset{\subset}{\lra}  \cdots \overset{\subset}{\lra} \DiscInf_{d,\leq \infty} = \DiscInf_d^{\oo}.\]
Restricting along these, we obtain a tower of categories
\[\ast \simeq \PSh^\un(\DiscInf_{d,\leq 0}^\oo) \longleftarrow \PSh^\un(\DiscInf_{d,\leq 0}^\oo) \longleftarrow \cdots \longleftarrow\PSh^\un(\DiscInf_{d,\leq \infty}^\oo) = \PSh^\un(\DiscInf_{d}^\oo).\]
From it we can recover $ \PSh^\un(\DiscInf_{d}^\oo)$ as the limit $\lim_{k \to \infty}  \PSh^\un(\DiscInf_{d,\leq k}^\oo)$. The terms in this limit can be understood inductively through the following results:
\begin{enumerate}[(a)]
	\item the description of the first stage as $\smash{\PSh^\un(\DiscInf_{d,\leq 1}^\oo)} \simeq \Spc_{/\BO(d)}$ and 
	\item the existence of pullback squares
\[\begin{tikzcd}\PSh^\un(\DiscInf_{d,\leq k}^\oo) \dar \rar &\Fun([2], \PSh(\BO(d) \wr \Sigma_k))\dar \\[-7pt]
	\PSh^\un(\DiscInf_{d,\leq k-1}^\oo) \rar &  \Fun([1],\PSh(\BO(d) \wr \Sigma_k))
\end{tikzcd} \qquad \text{for $k \geq 2$.}\]
\end{enumerate}
The latter is proven using a general recollement theorem, and informally say that to extend a presheaf on $\smash{\DiscInf_{d,\leq k-1}^\oo}$ to $\smash{\DiscInf_{d,\leq k}^\oo}$, we need to prescribe its value on exactly $k$ discs, compatibly with maps from a latching object and to a matching object and equivariantly for automorphisms of $\ul{k} \times \bfR^d$. This describes the layers of the tower as a category of fillers. Note that all descriptions of categories here are as limits and this, since mapping spaces in limit of categories are limits of mapping spaces, allows us to access maps between presheaves.

The naturality discussed in \cref{sec:structure} extends to these towers, and in fact we have already used the restriction to $\smash{\DiscInf_{d, \leq 1}^\oo}$ when discussing tangential structures.

\subsubsection{Manifold calculus and embedding calculus} Weiss first considered Disc-presheaves as a tool for studying manifolds, in his work on \emph{manifold calculus} and \emph{embedding calculus} \cite{WeissBulletin,Weiss,GoodwillieWeiss}. Rephrased in our setting, manifold calculus concerns the spaces 
\[X(M) \coloneqq \Map_{\PSh^\un(\DiscInf_d^\oo)}(E_M,X)\]
for $M \in \ManInf_d$ and $X \in \PSh^\un(\DiscInf_d^\oo)$, and embedding calculus the special case where $X = E_{N}$. Passing to mapping spaces, the cardinality tower gives for $M \in \ManInf_d$ and for $X \in \PSh^\un(\DiscInf_d^\oo)$ gives a tower conceptualised as approximations from which we recover $X(M)$ as the limit
\vspace{-.6cm}\begin{equation}\label{eqn:embcalc-tower-gen}
	\vspace{-.3cm}\begin{tikzcd}[row sep=0.2cm,/tikz/column 2/.append style={nodes={anchor=base west}}]
		&\hspace{2.52cm} \vdots \hspace{2.52cm} \ar[d, to path={-- (\tikztotarget.north -| \tikztostart)}]\\
		&T_2X(M) \coloneqq \Map_{\PSh^\un(\DiscInf_{d,\leq 2}^\oo)}(E_M,X)\ar[d, to path={-- (\tikztotarget.north -| \tikztostart)}]\\
		&T_1X(M) \coloneqq \Map_{\PSh^\un(\DiscInf_{d,\leq 2}^\oo)}(E_M,X)\ar[d, to path={-- (\tikztotarget.north -| \tikztostart)}]\\
		X(M) \coloneqq \Map_{\PSh^\un(\DiscInf_d^\oo)}(E_M,X) \rar\arrow[ur, bend left=10]\arrow[uur, bend left=17.5]\arrow[uuur, bend left=25,shift left=1]& T_0X(M) \coloneqq \Map_{\PSh^\un(\DiscInf_{d,\leq 0}^\oo)}(E_M,X) \simeq \ast.\\
	\end{tikzcd}
\end{equation}
Its initial stage and layers admit straightforward descriptions: firstly, we have
\[T_0X(M) \simeq \ast \qquad \text{and} \qquad T_1X(M) \simeq \Map_{/\BO(d)}(M,|X|)\]
where $|X| \in \Spc_{/\BO(d)}$ is the \emph{underlying space} of $X$ with vector bundle given by $|X| \coloneqq \smash{X(\bfR^d)_{\rm{O}(d)}} \to \smash{\BO(d)}$. Secondly, working out the recollement square one proves an equivalence
\[\fib_x\big(T_kX(M) \to T_{k-1}X(M) \big)\simeq \left\{\text{ fillers in} \begin{tikzcd} \partial^h C_k(M) \rar \dar & X(\ul{k} \times \bfR^d)_{O(d) \wr \Sigma_k} \dar \\[-7pt]
C_k(M) \rar \arrow[dashed]{ru} & {[\lim_{S \subsetneq \ul{k}} X(S \times \bfR^d)]}_{O(d) \wr \Sigma_k}\end{tikzcd} \right\}\]
where $C_k(M) \coloneqq \Emb^\oo(\ul{k},M)_{\Sigma_k} \subset M^k$ is the unordered configuration space and $\partial^h C_k(M)$ the homotopy link of the thick diagonal $\Delta_k(M) \coloneqq \{(x_1,\ldots,x_k) \mid x_i = x_j \text{ for some $i\neq j$}\} \subset M^k$.

\begin{ex}Let us take $X = \Emb^\oo(-,N)$ for a $k$-dimensional smooth manifold $N$, and abbreviate $\smash{T_k\Emb^\oo(M,N) \coloneqq \Map_{\PSh(\DiscInf_{d,\leq k}^\oo)}(E_M,E_N)}$ for $0 \leq k \leq \infty$. Then we get a tower
	\vspace{-.4cm}\begin{equation}\label{equ:embcalc-tower-mfd}
		\vspace{-.2cm}\begin{tikzcd}[row sep=0.2cm,/tikz/column 2/.append style={nodes={anchor=base west}}]
			&[10pt]\hspace{1cm}\vdots\hspace{1cm}\ar[d, to path={-- (\tikztotarget.north -| \tikztostart)}]\\
			&T_2\Emb^\oo(M,N)\ar[d, to path={-- (\tikztotarget.north -| \tikztostart)}]\\
			&T_1\Emb^\oo(M,N)\ar[d, to path={-- (\tikztotarget.north -| \tikztostart)}]\\
			T_\infty \Emb^\oo(M,N)\rar\arrow[ur, bend left=10]\arrow[uur, bend left=17.5]\arrow[uuur, bend left=25,shift left=1]&T_0\Emb^\oo(M,N) \simeq \ast.\\
		\end{tikzcd}
	\end{equation}
The first stage is $T_1 \Emb^\oo(M,N) \simeq \Map_{/\BO(d)}(M,N)$ with $M \to \BO(d)$ and $N \to \BO(d)$ classifying the tangent bundle: this is equivalent to the space of bundle isomorphisms $TM \to TN$ or by Smale--Hirsch theory, as long as $M$ has no closed components, the space of immersions $\rm{Imm}(M,N)$. The higher layers are given by spaces of fillers of a square involving only configuration spaces of $M$ and $N$; these are well-studied objects in homotopy theory.\end{ex} 

To elucidate why this is a ``calculus'', we exhibit these constructions as polynomial approximations as in \cite{Weiss}. An open cover $\scr{V}$ of $M$ is a \emph{$\scr{J}_k$-cover} if every collection of $\leq k$ points in $M$ is contained in some $V \in \scr{V}$ and it is \emph{complete} if it contains a $\scr{J}_k$-cover of every finite intersection $V_1 \cap \cdots \cap V_r$ for $V_1,\ldots,V_r \in \scr{V}$. Then the $T_kX(-)$ has \emph{descent for $\scr{J}_k$}-covers---equivalently, are sheaves for a Grothendieck topology of $\scr{J}_k$-covers---in the sense that the map induced by contravariant naturality of \eqref{eqn:embcalc-tower-gen} in $M$
\[T_kX(M) \lra \lim_{V \in \scr{V}} T_kX(V)\]
is an equivalence for each non-empty complete $\scr{J}_k$-cover $\scr{V}$ of a smooth $k$-dimensional manifold $M$. By \cite{BoavidadeBritoWeissSheaves}, having descent for $\scr{J}_k$-covers is equivalent to (i) taking increasing unions $U_1 \subset U_2 \subset \cdots \subset M$ of open subsets to limits (being \emph{good}), (ii) taking, for disjoint closed subsets $A_1,\ldots,A_{k+1} \subset M$ the cubical diagram $\ul{k+1} \supset S \mapsto M \setminus (\cup_{i \notin S} A_i)$ to a cartesian cubical diagram (being \emph{polynomial of degree $k$}); in fact, the $T_k X$ are the initial presheaves under $X$ that are good and polynomial of degree $k$.

\subsection{The second stage}\label{sec:second-stage} The first stage is well-understood, as $\PSh^\un(\DiscInf^\oo_{\leq 1}) \simeq \Spc_{/\BO(d)}$: a unital presheaf $X$ on $\smash{\DiscInf_{d,\leq 1}^\oo}$ (a \emph{1-truncated presheaf}) is the same as a space $|X|$ with a $d$-dimensional vector bundle $\xi$ over it; for $X = E_M$ this is $M$ with its tangent bundle. The study of the second stage goes back to work of Haefliger and Dax, see \cite{GoodwillieKleinWeissH}, and the additional data in a unital presheaf $X$ on $\smash{\DiscInf_{d,\leq 2}^\oo}$ (a \emph{2-truncated presheaf}) is the left $\Sigma_2$-equivariant commutative square
\[\begin{tikzcd} S(\xi) \rar \dar & X(\ul{2} \times \bfR^d)_{\rm{O}(d)^2} \dar \\[-7pt]
	{|X|} \rar{\Delta} & {|X| \times |X|}\end{tikzcd} \qquad \qquad \begin{tikzcd} S(TM) \rar \dar & M \times M \backslash \Delta \dar \\[-7pt]
	M \rar{\Delta} & M \times M\end{tikzcd}\]
where $S(\xi)$ is the $(d-1)$-dimensional spherical fibration associated to $\xi$. For $X = E_M$ this simplifies to the more familiar $\Sigma_2$-equivariant commutative square that appears on the right. This is intimately related to three properties shared by Disc-presheaves $E_M$ arising from a closed manifold $M$ (there are also versions for compact manifolds with boundary), which appear in recent work of Naef and Safronov, and of myself, Krannich, and Mezher:
\begin{enumerate}[\noindent (1)]
	\item \emph{Poincar\'e duality.} \label{enum:poincare} One can construct from a square of the left form a zigzag of maps of spectra from $\bfS$ to $\Th(-\xi)$, and if $|X|$ is a Poincar\'e duality space one can ask (i) whether the wrong-way map that appears in it is an equivalence and (ii) whether the resulting map $\bfS \to \Th(-\xi)$ exhibits the stabilisation of $S(\xi)$ as the Spivak normal fibration. Both are true for $X = E_M$ if $M$ is a closed manifold \cite{PriggeThesis,NaefSafronov,KKM}.
	\item \emph{Finiteness structure.} \label{enum:finiteness} If $|X|$ is a compact object in $\Spc$ then has a canonical class in Waldhausen's algebraic $K$-theory of spaces $\rm{A}(X)$ and a finiteness structure on this is a lift along the assembly map $\rm{A}(*) \otimes X \to \rm{A}(X)$. There is a trace map $\rm{tr} \colon \rm{A}(X) \to \rm{THH}(X)$ and a 2-truncated presheaf with Poincar\'e duality has a preferred lift of its image along the $\rm{THH}$-assembly map $\bfS \otimes X \to \rm{THH}(X)$ \cite{NaefSafronov}. 
	\item \emph{Stratification of products.} \label{enum:stratification} The left square yields a map $|X| \cup_{S(\xi)} \smash{X(\ul{2} \times \bfR^d)_{\rm{O}(d)^2}} \to |X| \times |X|$. For the right square, i.e.~for a manifold $M$, this is an equivalence and it corresponds to the natural stratification of the product $M^2$ with strata the ordered configuration spaces.
\end{enumerate} 
\cref{enum:stratification} admits a generalisation for $k$-truncated presheaves: the left Kan extension of $E_M$ along the functor $\pi_0 \colon \DiscInf_d^{\oo} \to \Fin$ satisfies a Segal condition. In the more algebraic setting of \cite{WillwacherObstruction}, this is referred to as \emph{being of configuration space type}. It remains open what the correct generalisations of (1) and (2) are. The underlying problem of these considerations is to characterise those Disc-presheaves that arise from manifolds:

\begin{prob}Determine the essential image of the map \eqref{eqn:emap} and its truncated variants in terms internal to the category of presheaves, and relate it to classical invariants of manifolds.
\end{prob}

An alternative perspective on this problem is that it should tell us which Disc-presheaves should be considered \emph{as good as} manifolds. We should then classify these, or use custom-built ones, not arising from manifolds, in our arguments.

\subsection{Outlook: what is next for the tower?}

\subsubsection{Develop more calculational tools} \label{sec:computations} Towers like \eqref{eqn:embcalc-tower-gen} and \eqref{equ:embcalc-tower-mfd} have associated Bousfield--Kan spectral sequences. These are particularly useful for establishing qualitative properties, as is done in some of the applications discussed in \cref{sec:applications}. More quantitative computations, however, have only been possible under restricted circumstances: in a range of degrees, for particularly domains and targets, or localised or completed.

Particular progress has been made for rationalised embedding calculus; I will give a purposely vague description so as to cover a variety of results, such as \cite{ALTV,AroneTurchin,FTWRn,WillwacherObstruction}. The idea is to combine the operadic perspective of \cref{sec:operadic} with Kontsevich's formality for the little $d$-discs operad $E_d$ \cite{KontsevichFormality}, to find $L_\infty$-algebras of a graphical nature whose Maurer--Cartan spaces compute mapping spaces of rationalised presheaves:
\[\rm{MC}(\rm{Graph}_{M,N}) \simeq \Map_{\PSh(\DiscInf_d^\oo)}(E_M^\bfQ,E_N^\bfQ).\]
There are then separate problems to actually compute with these graphical models, and to relate the right side to the rationalisations of the original mapping spaces.

\begin{rem}Configuration space integrals \cite{KontsevichFeynmann}, as used recently by Watanabe (see e.g.~\cite{WatanabeS4}), are invariants of manifold bundles that are indexed by the homology of certain graph complexes. These are surely closely related to embedding calculus, though no comparison has appeared as of yet.\end{rem}

We do not have similarly powerful descriptions when working with other localisations or completions, or working stably.

\begin{question}What are integral or spectral analogues of these graphical models?
\end{question} 

\subsubsection{The tower versus its limit} \label{sec:tower-vs-limit} The cardinality tower played a crucial conceptual role in the development of embedding calculus and its applications. It is, however, not clear whether we should take this as a sign to make it a more central object of study, or look past it. The correct answer likely depends on the situation.

\smallskip

Favouring the tower over its limit is that taking the limit of the tower is a main source of difficulty, as this does not interact well with many constructions and properties: this is the reason for truncations appearing in the results of \cite{Mezher,KKM} (see \cref{thm:smooth} \eqref{enum:smooth-2}). To sidestep this, one could consider the tower itself as a pro-space
\[\{\Map_{\PSh(\DiscInf_{d,\leq k}^\oo)}(X,Y)\}_{k \geq 1} \in \Pro(\Spc).\]
The computations using the tower can then be interpreted as determining this pro-object, and it becomes a separate problem to take its limit. One advantage is that it interacts better with localisations: e.g., the limit of the rationalisation of the pro-object rather than the rationalisation of its limit, is studied in rationalised embedding calculus as in \cref{sec:computations}.

\smallskip

There several reasons to favour the limit over the tower. Firstly, Disc-presheaves that are not truncated have properties, such as the isotopy extension result of \cref{sec:isotopy-extension}, that are not shared by truncated Disc-presheaves. This makes it plausible they may be studied by techniques---possibly of a more surgery-theoretic nature---that can not have truncated analogues. Secondly, rational computations using graphical models suggest that the tower is in fact a particularly bad filtration from the point of view of computations, since it involves extreme amounts of cancellation. It also does not interact well with additional structure like the loop-order filtration, though \cite{WillwacherTruncated} suggests that as a pro-object it is better behaved. At the very least, it is worthwhile to investigate other approaches, e.g.~alternative towers or the incorporation of other forms of calculus---such as Weiss' orthogonal calculus c.f.~\cite{AroneDerivatives}.

\section{Desideratum \eqref{enum:des-sdisc}: understanding the difference between manifolds and Disc-presheaves} \label{sec:sdisc} Having the tools to understand Disc-presheaves and relate them to manifolds, we now turn to understand the difference.

\subsection{Convergence}\label{sec:convergence} When using Disc-presheaves to study manifolds, we are at first most interested in those Disc-presheaves associated to manifolds and hence whether the functor $E \colon \ManInf_d^\oo \to \PSh(\DiscInf_d^\oo)$ is fully faithful. That is, is the left comparison map below an equivalence?
\vspace{-.2cm}\[\Emb^\oo(M,N) \xrightarrow{\rm{comp}} T_\infty \Emb^\oo(M,N) \overset{\simeq}\lra \lim_{k \to \infty} T_k \Emb^\oo(M,N).\]

\subsubsection{Convergence results} The foundational result in this direction was obtained by Goodwillie and Weiss \cite{GoodwillieWeiss}, building on work of Goodwillie and Klein \cite{Goodwillie,GoodwillieKleinI,GoodwillieKleinII}, and in high and low dimensions improved by myself and Krannich \cite{KKSurfaces, KKOperadic}:

\begin{thm}\label{thm:convergence} Suppose that $M$ and $N$ are compact $d$-dimensional smooth manifolds, possibly with boundary. Then the comparison map $\rm{comp}$ is an equivalence if
	\begin{enumerate}
		\item \label{enum:convergence-1} the handle dimension of $M$ is $\leq d-3$, or
		\item \label{enum:convergence-3} the dimension is $d \leq 2$, or
		\item \label{enum:convergence-2} the inclusion $\partial M \to M$ is an equivalence on tangential 2-types and $d \geq 5$.
	\end{enumerate}
\end{thm}

The statement needs some explanation: $M$ has \emph{handle dimension $\leq h$} if it can be constructed by iterated handle attachments of index $\leq k$, and $\partial M \to M$ is an \emph{equivalence on tangential 2-types} if it induces an equivalence on fundamental groupoids and for all components the restricted second Stiefel--Whitney class $w_2 \colon \pi_2(\partial M) \to \bfZ/2$ is non-trivial if and only if $w_2 \colon \pi_2(M) \to \bfZ/2$ is non-trivial. Note that (1) implies (3) when $d \geq 5$. There are variants for triads and embeddings relative to part of the boundary, and for topological embeddings between smoothable topological manifolds (using the smoothing theory of \cref{sec:smoothing-theory}).

If the conclusion of \cref{thm:convergence} holds, we say that embedding calculus \emph{converges}. Convergence results can be used in two directions. On the one hand, they give a homotopy-theoretic description of certain spaces of embeddings, accessible through the results in \cref{sec:discpresheaves}. On the other hand, they give us access to mapping spaces between certain Disc-presheaves. 

\begin{rem}Goodwillie and Weiss actually proved a stronger result, namely that $\Emb^\oo(-,N)$ is \emph{analytic}. This in particular implies that under hypothesis \eqref{enum:convergence-1} of \cref{thm:convergence}, the map $\Emb^\oo(M,N) \to T_k \Emb^\oo(M,N)$ increases in connectivity with $k$. It seems better to separate convergence per se from the rate of convergence.
\end{rem}

\subsubsection{Revisiting convergence} \label{sec:convergence-proof} \cref{thm:convergence} describes the state of convergence results as of writing, and it leaves something to be desired. Firstly, we do not know convergence for topological embeddings when the target is not smoothable. Arguably, the lack of this technical tool is the main distinction between topological manifold theory and smooth manifold theory. Secondly, the current proof is a ``first generation'' one and rather complicated, relying on hard multiple disjunction results. There has been no significant improvements or simplifications to the arguments of \cite{Goodwillie,GoodwillieWeiss,GoodwillieKleinI,GoodwillieKleinII}, nor any applications of similar techniques in other settings that have shed more light on them. 

\begin{prob}Find a ``second generation'' proof of convergence of embedding calculus, which also applies to the piecewise-linear and topological settings. 
\end{prob}


\subsection{Isotopy extension}\label{sec:isotopy-extension} \cref{thm:convergence} does not apply to closed manifolds and their diffeomorphisms. However, it is possible to transfer information about spaces of embeddings, by comparing them using isotopy extension in the shape of a fibre sequence
\[\Emb^\oo_\partial(M-\rm{int}(P),N) \lra \Emb^\oo(M,N) \lra \Emb^\oo(P,N)\]
for $P \subset M$ a compact codimension zero submanifold and fibre taken over the inclusion. There is an analogue of this for maps of presheaves \cite{KKnudsen}: 

\begin{thm}\label{thm:isotopy-extension} If the map $\Emb^\oo(\ul{k} \times \bfR^d \sqcup P,N) \to T_\infty \Emb^\oo(\ul{k} \times \bfR^d \sqcup P,N)$ is an equivalence for all $k \geq 0$, then there is a fibre sequence
	\[T_\infty \Emb^\oo_\partial(M-\rm{int}(P),N) \lra T_\infty \Emb^\oo(M,N) \lra T_\infty \Emb^\oo(P,N).\]
\end{thm}

This is a powerful tool: e.g., it is used to reduce \cref{thm:convergence} \eqref{enum:convergence-2} to \eqref{enum:convergence-1} by handle manipulation, and it is used to show that $E$ is fully-faithful in low dimensions:

\begin{ex}\label{ex:low-dimensions} The Cerf--Gramain approach to diffeomorphism groups of compact surfaces determines these using only the constructions of \cref{sec:structure} and isotopy extension \cite{Gramain}. They may thus be replicated in the category of Disc-presheaves to yield the case $d=2$ of \cref{thm:convergence} \eqref{enum:convergence-3} \cite{KKSurfaces}; the cases $d=0,1$ are similar but easier. Can Disc-presheaves be used to study mapping class groups, e.g.~the Johnson filtration or Torelli Lie algebra?
\end{ex}

\subsection{Non-triviality of Disc-structure spaces} \label{sec:non-triviality} Recall from \cref{sec:moduli-sdisc} that for closed manifolds, or more generally for compact manifolds with given closed $(d-1)$-dimensional smooth manifold $P$ as their boundary, the difference between moduli spaces of manifolds and moduli spaces of Disc-presheaves is measured by the Disc-structure spaces
\[S^\Disc_\partial(X) \coloneqq \fib_X\big[\ManInf_{d}^P \overset{E}\lra \rm{LMod}_{E_{P \times \bfR}}(\PSh(\DiscInf^\oo_d))^\simeq\big],\]
where $\ManInf_{d}^P$ denotes the groupoid whose objects are compact smooth manifolds $M$ whose boundary is identified with $P$ and whose mapping spaces are spaces of diffeomorphisms relative to $P$. \cref{thm:convergence} does not apply here unless $d\leq 2$, in which case the Disc-structure spaces are thus contractible. In fact, by work of Krannich and myself these are \emph{only} dimensions in which the Disc-structure spaces are always contractible \cite{KKSDisc}:

\begin{thm}\label{thm:non-triviality} Suppose $d \neq 4$, then $S^\Disc_\partial(D^d)$ is contractible if and only if $d \leq 2$.
\end{thm}

Let us outline the argument for large dimensions, based on \cite{KKSquare}. By smoothing theory---classically in the form of Morlet's theorem $\smash{\BDiff_\partial(D^d) \simeq \Omega^d_0 \tfrac{\Top(d)}{\rm{O}(d)}}$ and for Disc-presheaves in the form of \eqref{eqn:delooping}---we get an equivalence
\[S^\Disc_\partial(D^d) \simeq \fib\big[\Omega^d_0 \tfrac{\Top(d)}{\rm{O}(d)} \to \Omega^d_0 \tfrac{\Aut(E_d)}{\rm{O}(d)}] \simeq \Omega^{d+1} \tfrac{\Aut(E_d)}{\rm{Top}(d)},\]
casting the problem in terms of the map $\BTop(d) \to \BAut(E_d)$ encoding the action of $\Top(d)$ on the little $d$-discs operad $E_d$. The idea is that this map must be far from an equivalence. There are two possible worlds: (a) $\BAut(E_d)$ has an uncountable homotopy group, (b) it does not. In case (a), this map can not be an equivalence because all homotopy groups of $\BTop(d)$ are countable. In case (b), one can deduce there is an equivalence $\BSAut(E_d)_\bfQ \simeq \smash{\BSAut(E^\bfQ_d)}$ between the rationalisation of the automorphisms and automorphisms of the rationalisation. We then use that increasing the dimension has a different effect for homeomorphisms than for automorphisms of the rationalised $E_d$-operad: there is a commutative square
\[\begin{tikzcd} \BSTop(d-2) \rar \dar & \BSAut(E_{d-2}^\bfQ) \dar{\simeq \ast} \\[-7pt]
	\BSTop(d) \rar & \BSAut(E_d^\bfQ) \end{tikzcd} \]
where the left vertical map is highly-connected but the right vertical map can be shown to be null-homotopic, making it impossible that the bottom map is a rational equivalence. Though this is a proof, it is unsatisfying that we have to argue by cases:

\begin{question}Are any of the homotopy groups of $\BAut(E_d)$ uncountable?\end{question}

This is related to the question whether the pro-space $\{\BSAut(E_{d,\leq k})\}_{k \geq 1}$ of automorphisms of the truncated little $d$-discs operad has pro-constant or Mittag--Leffler homotopy groups. 

\subsubsection{Outlook: relation to orthogonal calculus}
One may wonder whether---in analogy with \cite[Section 3.3]{RandalWilliamsDiscs} and the perspective of \cite{KKOdd}---there is a pullback square
\[\begin{tikzcd} \BTop(d) = \mathbf{Bt}(\bfR^d) \rar \dar & P_1 \mathbf{Bt}(\bfR^d) \dar \\[-7pt]
	\BAut(E_d) = \mathbf{Ba}(\bfR^d) \rar & P_1 \mathbf{Ba}(\bfR^d), \end{tikzcd}\]
in terms of the orthogonal functors $\mathbf{Bt} \coloneqq (V \mapsto \BTop(V))$ and $\mathbf{Ba} \coloneqq (V \mapsto \BAut(E_V))$ (not rationalised, unlike in loc.cit.) and the first stage $P_1$ of the orthogonal calculus tower, so that the only source of any such uncountable homotopy groups are the zeroth and first orthogonal derivatives. A starting point would be to understand the map on the zeroth stages:
\[P_0 \mathbf{Bt}(\bfR^d) = \underset{d \to \infty} \colim\, \BTop(d) \lra P_0 \mathbf{Ba}(\bfR^d) = \underset{d \to \infty} \colim\, \BAut(E_d).\]

\begin{prob} Determine $\underset{d \to \infty} \colim\, \BAut(E_d)$.
\end{prob}

\begin{rem}Randal-Williams has suggested to combine this with \cref{sec:tower-vs-limit} and consider the $\Pro(\Spc)$-valued orthogonal functor $\mathbf{Ba}^\rm{pro} \coloneqq (V \mapsto \{\BAut(E_{V, \leq k})\}_{k \geq 1})$.\end{rem}

\begin{ex}The topological Hirzebruch $L$-classes in $H^*(\BTop;\bfQ)$ are not pulled back from $H^*(\colim_{d \to \infty}\,\BAut(E_{d,\leq k});\bfQ)$ for fixed finite $k$---the family signature theorem fails in truncated embedding calculus \cite{PriggeThesis}. However, it is possible that they are pulled back from $H^*(\colim_{d \to \infty}\,\BAut(E_d);\bfQ)$ and there are detected on a copy of the $p$-adic integers in the homotopy groups arising as a limit of finite abelian groups in the truncations.\end{ex}

\subsection{Structural properties of Disc-structure spaces} As the Disc-structure spaces are often non-trivial, we should find properties that will help us understand them. Combining isotopy extension for spaces of embeddings and maps of Disc-presheaves (i.e.~\cref{thm:isotopy-extension}), Krannich and myself related the Disc-structure spaces of different manifolds, obtaining an invariance result \cite{KKSDisc}, similar to that for spaces of positive scalar curvature metrics \cite{ERW,EbertWiemeler}:

\begin{thm}\label{thm:2-type-invariance} Let $M$ and $N$ be compact $d$-dimensional smooth manifolds with $d \geq 5$. If they have equivalent tangential 2-types, then there is an equivalence
	\[S^\Disc_\partial(M) \simeq S^\Disc_\partial(N).\]
\end{thm}

What it means for two manifolds to have same tangential 2-type was discussed following \cref{thm:convergence}, and it has two implications. Firstly, if $M$ and $N$ are homeomorphic they have the same tangential 2-type so Disc-structure spaces only depend on the underlying topological manifold; this can also be proven directly using the smoothing theory of \cref{sec:smoothing-theory} \cite{KKOperadic}. Secondly, it implies that taking a connected sum with $S^k \times S^{d-k}$ for $2 \leq k \leq d-2$ does not affect Disc-structure spaces. Combining this with the machinery of operads with homological stability \cite{BBPTY}, one can deduce a second structural property \cite{KKSDisc}:

\begin{thm}\label{thm:infinite-loop-space} Let $M$ be a compact $d$-smooth manifold with $d \geq 8$. Then $S^\Disc_\partial(M)$ admits the structure of an infinite loop space.
\end{thm}

From the perspective of \cref{sec:non-triviality}, this begs the following question:

\begin{question}Does $\tfrac{\rm{Aut}(E_d)}{\Top(d)}$ admit the structure of an infinite loop space?
\end{question}

\subsection{Outlook: a higher-algebraic description of Disc-structure spaces} The structural properties of the Disc-structure spaces are strongly reminiscent of two objects that appear in the classical approach to manifolds using surgery theory and pseudoisotopy theory (c.f.~the first two points on the list in \cref{sec:second-stage}):
\begin{enumerate}[(1)]
	\item \emph{L-theory}. The main result of surgery theory can be roughly phrased as a the existence of a fibre sequence
	\[S(X) \lra N(X) \lra L(X),\]
	involving a structure space $S(X)$ of manifold structures on a Poincar\'e complex, a space of normal invariants $N(X)$, and a $L$-theory space $L(X)$. The latter is the infinite loop space of an $L$-theory spectrum, which by the $\pi$-$\pi$-theorem only depends on the fundamental groupoid of $X$ and the first Stiefel--Whitney class \cite{WallSCM}.
	\item \emph{The fibre of the cyclotomic trace}: The main result of pseudoisotopy theory can be phrased as the existence of a highly-connected map
	\[\scr{H}(M) \lra \Omega^{\infty+1} \rm{Wh}^\Diff(M)\]
	from a space of stable $h$-cobordisms starting at a manifold to the infinite loop space of the Whitehead spectrum of $M$, a summand of the algebraic $K$-theory of spaces spectrum $\rm{A}(M)$. The latter is nowadays often studied using trace methods, in particular through a cyclotomic trace $\rm{trc} \colon \rm{A}(M) \to \rm{TC}(M)$ to topological cyclic homology. Its fibre only depends on the fundamental groupoid of $M$ \cite{DGM}.
\end{enumerate}
A relationship between Disc-structure spaces and $L$-theory appeared in work of Krannich, myself, and Mezher, as follows: the normal invariant construction associated to the second stage, described in \cref{sec:second-stage} \eqref{enum:poincare}, should induce for a compact smooth manifold $M$ a map only depending on the tangential 2-type
\[S^\Disc_\partial(M) \lra \Omega L(M).\]

\begin{question}How non-trivial is this map?\end{question}

The relationship of Disc-structure spaces to topological cyclic homology is work-in-progress by Krannich and Mu\~{n}oz-Ech\'{a}niz, and is related to the work of Naef and Safronov (described in \cref{sec:second-stage} \eqref{enum:finiteness}). It is worth stating the implicit problem, necessarily vague:

\begin{prob}\label{prob:higher-algebra} Give a description of Disc-structure spaces in terms of ``higher-algebraic'' constructions with ``good'' computational properties.
\end{prob}

An answer for this in terms of known objects would be very powerful. However, it would even be better if the answer involves novel objects, e.g.~new variants of algebraic $K$-theory, with analogues of trace maps that involve spaces of graphs by a suggestion of Goodwillie. As stated in the introduction, whether \cref{prob:higher-algebra} admits an interesting answer should be considered a measure of the success of Disc-presheaves.

\section{Applications}\label{sec:applications} In addition to inspiring developments in other fields, the measure of success of Disc-presheaves as an approach to the study of manifolds should be whether one can actually use them to study manifolds. It should be clear that the answer is \emph{yes}, if only by the convergence results of \cref{thm:convergence}. To support this we give several recent examples of breakthroughs that rely on Disc-presheaves (or embedding calculus). Many other results had to be left out, including computations for long knots, e.g.~\cite{FTW,BoavidadeBritoHorel}, barbell diffeomorphisms and their applications, e.g.~\cite{BudneyGabai}, the Alexander trick for contractible manifolds \cite{GRWAlexander,KKOperadic},  and the connection to Whitehead torsion and string topology \cite{NaefSafronov}.

\subsection{Advantages relative to surgery theory and pseudoisotopy theory} Before going into specifics, let us explain on a high level why using Disc-presheaves has applications that were not previously accessible through the classical method of combining \emph{surgery theory} and \emph{pseudoisotopy theory} (c.f.~\cref{sec:second-stage}). This approach, which compares manifolds to block manifolds and these in turn to Poincar\'e complexes, is of course highly successful: it provides much insight in the classification of high-dimensional manifolds and inspired major developments in algebraic $K$- and $L$-theory. However, it has two drawbacks:
\begin{enumerate}[(a)]
	\item it only provides information about a truncation in a range that depends on the dimension,
	\item it involves two steps leading to additional extension problems (though see \cite{WeissWilliams}).
\end{enumerate}
Using Disc-presheaves avoids these drawbacks. Of course, an approach using Disc-presheaves is unlikely to replace the classical one, and the stronger results can be often obtained by combining both.

\subsection{Distinguishing homotopy spheres} A first test of the strategy of studying manifolds through their Disc-presheaves is whether they can distinguish homotopy $d$-spheres $\Sigma$. This is harder than expected, since the 1-truncations $\smash{E_{S^d}, E_{\Sigma} \in \PSh(\DiscInf_{d,\le 1}^\oo)}$ are equivalent by the classification of tangent bundles to homotopy spheres; in fact, each of the individual spaces $\Emb(\ul{k} \times \bfR^d,\Sigma)$ is equivalent to $\Emb(\ul{k} \times \bfR^d,S^d)$. However, the higher truncations can distinguish some of them: the first examples are a consequence of convergence \cite{KKnudsen}, using examples of homotopy spheres that do not admit an embedding of low codimension into an Euclidean space \cite{HLS}, unlike $S^d$ which embeds in $\bfR^{d+1}$.

For the more comprehensive answer given in \cite{KKM}, recall the Kervaire--Milnor classification describes the group $\Theta_d$ of oriented homotopy $d$-spheres, outside finitely many exceptional dimensions, via an short exact sequence
\[0 \lra \bP_{d+1} \lra \Theta_d \overset{\rm{PT}}\lra \coker(J)_d \lra 0,\]
in terms of its subgroup $\bP_{d+1}$ of those homotopy $d$-spheres that bound a stably parallelisable compact smooth manifold, which is a cyclic finite group of known order, and the cokernel of the stable $J$-homomorphism. The latter can be rephrased in terms of the construction in \cref{sec:second-stage} \eqref{enum:poincare} yielding part \eqref{enum:smooth-1} of the following result (which is then also true for $k = \infty$):

\begin{thm}\label{thm:smooth} Let  $\Sigma,\Sigma' \in \Theta_d$ and $1<k<\infty$.
	\begin{enumerate}[(i)]
		\item \label{enum:smooth-1}  If $\rm{PT}(\Sigma)-\rm{PT}(\Sigma') \neq 0 \in \coker(J)_d$, then $E_{\Sigma} \not \simeq E_{\Sigma'} \in \PSh(\DiscInf_{d,\leq k}^\oo)$.
		\item \label{enum:smooth-2}  If $\Sigma-\Sigma' \in \bP_{d+1}$ and $d {\not \equiv} 1 \text{ (mod 4)}$, then $E_{\Sigma} \simeq E_{\Sigma'} \in \PSh(\DiscInf_{d,\leq k}^\oo)$.
	\end{enumerate}
\end{thm}

\noindent The second part uses (a) that any $(4k-1)$-sphere $\Sigma$ can be built by gluing two solid handlebodies $\natural_g D^{2k} \times S^{2k-1}$ along $W_g \coloneqq \sharp_g S^{2k-1} \times S^{2k-1}$ using a diffeomorphism that lies in the finite residual of the oriented mapping class group $\pi_0\,\Diff^+(W_g)$ of $W_g$, (b) a computation of this finite residual building on work of Krannich and Randal-Williams \cite{Krannichmcg,KRWmcg}, and (c) a surprising result of Mezher that the group $\pi_0 \, T_k \Diff^+(W_g)$ is residually finite \cite{Mezher}. It remains open whether the analogue of \cref{thm:smooth} \eqref{enum:smooth-2}  holds if $d {\equiv} 1 \text{ (mod 4)}$ or if $k=\infty$.

\subsection{The Weiss fibre sequences and its first applications} There is a class of applications that is based on a fibre sequence relating diffeomorphism groups and embedding spaces, originally due to Weiss in a special case \cite{WeissPontryagin} (see \cite{Cardim} for a precursor). 

\subsubsection{The Weiss fibre sequence} Suppose we are given a manifold triad $M$, its boundary subdivided into two codimension zero submanifolds $\partial_0 M$ and $\partial_1 M$ meeting at corners. Then, as in \cref{sec:disjoint-union-gluing}, we can think of a collar on $\partial_1 M$ as an $E_1$-algebra in $\ManInf_d$ and $M$ as a left $E_1$-module over it. On the level of moduli spaces, taking the quotient by this action yields the \emph{Weiss fibre sequence}
\[\BDiff_\partial(M) \lra \BEmb^{\oo,\cong}_{\partial_0}(M,M) \lra \rm{B}^2\Diff_\partial(\partial_1 M \times [0,1]),\]
where the superscript $\cong$ indicates we only take those self-embeddings isotopic to a diffeomorphism fixing the boundary pointwise, and in the right term the second delooping is with respect to stacking in the interval direction \cite{KupersThesis}. 

Now that we have made a space of smooth embeddings appear, we can try to apply embedding calculus to it. For example, if we suppose that $M$ is 1-connected spin of dimension $d \geq 5$ and $\partial_1 M = D^{d-1}$, then a variant of the convergence result \eqref{enum:convergence-3} of \cref{thm:convergence} for embeddings of triads applies and we can identify the middle term of the Weiss fibre sequence with (certain components of) a mapping space of Disc-presheaves: the result is a fibre sequence
\begin{equation}\label{equ:embcalc-weiss-fibre-sequence} \BDiff_\partial(M) \lra \rm{B}T_\infty \Emb^{\oo,\cong}_{\partial_0}(M,M) \lra \rm{B}^2\Diff_\partial(D^d)\end{equation}
with $\rm{B}T_\infty \Emb^{\oo,\cong}_{\partial_0}(M,M)\coloneqq \smash{\rm{B}\Map^{\cong}_{\rm{LMod}_{E_{\partial_0 M \times \bfR}}(\PSh(\DiscInf_d^\oo))}(E_M,E_M)}$.

\subsubsection{Finite generation results}  \label{sec:finiteness} We explain a prototypical application of the Weiss fibre sequence \eqref{equ:embcalc-weiss-fibre-sequence}, which establishes that  the homotopy groups of certain diffeomorphism groups are finitely-generated in each degree.

It turns out to be a good choice to take $M$ to be $W_{g,1} \coloneqq D^{2n} \sharp (S^n \times S^n)^{\# g}$ for $2n \geq 6$, with the intent of letting $g \to \infty$. In that case, the left term $\BDiff_\partial(W_{g,1})$ of the fibre sequence
\begin{equation}\label{eqn:weiss-fiber-sequence-wg} \BDiff_\partial(W_{g,1}) \lra \rm{B}T_\infty \Emb^{\oo,\cong}_{\partial_0}(W_{g,1},W_{g,1}) \lra \rm{B}^2\Diff_\partial(D^{2n})\end{equation}
is homologically accessible in a range of degrees tending to infinity with $g$ through the work of Galatius and Randal-Williams on homological stability and stable moduli spaces \cite{GRWStable,GRWStabI,GRWStabII}. If we use embedding calculus to access the middle term, we have information about two of the three terms in \eqref{eqn:weiss-fiber-sequence-wg} and thus also about the remaining term $\rm{B}^2\Diff_\partial(D^{2n})$. Taking loops once, we learn about $\BDiff_\partial(D^{2n})$ in a range of degrees tending to infinity with $g$, which is \emph{no} limitation at all since it is independent of $g$. 

It is not hard to use the above strategy to establish that the homotopy groups of $\smash{\BDiff_\partial(D^{2n})}$ are finitely-generated in each degree when $2n \geq 6$. Similarly using odd-dimensional results on homological stability and stable moduli spaces due to Botvinnik and Perlmutter \cite{BotvinnikPerlmutter,Perlmutter} to get the cases $d=2n+1 \geq 9$, this yields \cite{KupersThesis}:

\begin{thm}For $d \neq 4,5,7$ the homotopy groups of $\BDiff_\partial(D^d)$ are finitely-generated in each degree.
\end{thm} 

This result can be fed back into \eqref{equ:embcalc-weiss-fibre-sequence} and combined with embedding calculus, to prove finite generation results for the homotopy groups of the left term $\BDiff_\partial(M)$ for certain smooth manifolds $M$; this approach was elaborated upon in joint work of Bustamante, Krannich, and myself \cite{BKK}. Future work of Krannich and myself will extend this strategy to other properties than finite-generation, via a suitable theory of Serre classes.

\subsubsection{Rational homotopy groups of diffeomorphisms of discs}
Of course, we should use the Weiss fibre sequence for more than merely qualitative results and perform some explicit computations. As explained well in \cite{RandalWilliamsDiscs}, the series of papers \cite{KRWTorelli,KRWAlgebraic,KRWKoszul,KRWDiscs} of myself and Randal-Williams does so using \eqref{eqn:weiss-fiber-sequence-wg} and uses it to determine the rational homotopy of diffeomorphism groups of even-dimensional discs outside certain bands. This essay does not provide enough space to go into the details.

\subsubsection{Further applications of the Weiss fibre sequence} The Weiss fibre sequence also underlies the work of Bustamante and Randal-Williams on diffeomorphisms of solid tori $S^1 \times D^{2n-1}$ \cite{BRWTori}, the work of Krannich on a homological approach to pseudoisotopy theory in the 1-connected case \cite{KrannichHomological}, and the work of Krannich and Randal-Williams on the rational homotopy of diffeomorphism groups of odd-dimensional discs \cite{KKOdd}; in these, the role of Disc-presheaves is taken on by the multiple disjunction results of Goodwillie--Klein that are essentially equivalent to the convergence result of \cref{thm:convergence} \eqref{enum:convergence-1}.

\subsection{Pontryagin--Weiss classes and a pullback decomposition} It is also possible to use the Weiss fibre sequence in the topological setting, with striking applications to characteristic classes of topological $\bfR^d$-bundles; this was Weiss' original motivation.

\subsubsection{The Weiss fibre sequence for homeomorphisms} If we consider the Weiss fibre sequence \eqref{equ:embcalc-weiss-fibre-sequence} for \emph{homeomorphisms} of triads, the right term of the Weiss fibre sequence is contractible by the Alexander trick and we get an equivalence
\begin{equation}\label{eqn:bhomeo-as-embcalc} \BHomeo_\partial(M) \overset{\simeq}\lra \rm{B}T_\infty \Emb^{\ot,\cong}_{\partial_0}(M,M)\coloneqq\rm{B}\Map^{\cong}_{\rm{LMod}_{E^\ot_{\partial_0 M \times \bfR}}(\PSh(\DiscInf_d^\ot))}(E^\ot_M,E^\ot_M)\end{equation}
between the classifying spaces of a homeomorphism group and of (certain components of) a mapping space of topological Disc-presheaves. This should be surprising, as it is a purely homotopy-theoretic descriptions---involving only tangential data and configuration spaces---of certain moduli spaces of topological manifolds.

\subsubsection{Pontryagin--Weiss classes} \label{sec:pontryagin-weiss} Weiss had the ingenious idea to use this to study rational topological Pontryagin classes \cite{WeissPontryagin}. It follows from smoothing theory and the finiteness of the groups $\Theta_d$ of homotopy spheres that the map induced by the inclusions
\[\BO \coloneqq \underset{d \to \infty}\colim\,\BO(d) \lra \BTop \coloneqq \underset{d \to \infty}\colim\,\BTop(d)\]
is a rational equivalence. In particular, the rational cohomology $H^*(\BTop;\bfQ)$ is given by a polynomial ring $\bfQ[p_1,p_2,\ldots]$ in the rational Pontryagin classes. It is well-known how Pontryagin classes behave under restriction to $\BO(d)$: namely, $p_k$ vanishes if and only if $k > d/2$. Weiss showed this is \emph{not} true for restriction to $\BTop(d)$, by constructing examples of topological $\bfR^d$-bundles over $S^{4k}$ satisfying $k > d$ whose $k$th Pontryagin class is non-zero. That is, not only do we have $p_k \neq 0 \in H^{4k}(\BTop(d);\bfQ)$ but this is detected on the image of the Hurewicz map $\pi_{4k}(\BTop(d)) \to H_{4k}(\BTop(d))$. The first part of this statement was generalised vastly in work of Galatius and Randal-Williams \cite{GRWPontryagin}: they proved that all $p_i$ as well as the Euler class $e$ are algebraically independent in $H^*(\BTop(d);\bfQ)$ as long as $d \geq 6$:
\[\bfQ[e,p_1,p_2,\cdots] \subseteq H^*(\BTop(d);\bfQ) \qquad \text{for $d \geq 6$.}\]
Their argument uses finite group actions and yields no information about the second part, i.e.~whether one can find bundles with non-zero Pontryagin classes \emph{over spheres}. Generalising Weiss' approach, Krannich and I answered this \cite{KKSquare}:

\begin{thm}\label{thm:pontryagin-weiss} For $d \geq 6$ and $k \geq 1$ there exists a topological $\bfR^d$-bundle over $S^{4k}$ whose $k$th Pontryagin class is non-zero.
\end{thm}

Equivalently, the map $\BTop(d) \to \BTop$ is surjective on rational homotopy groups for $d \geq 6$. Through Morlet's theorem $\BDiff_\partial(D^d) \simeq \Omega^d_0 \smash{\tfrac{\Top(d)}{\rm{O}(d)}}$ this gives the first known infinite families of non-zero rational homotopy groups of diffeomorphism groups of discs.

\subsubsection{A pullback decomposition} The proof of \cref{thm:pontryagin-weiss} contains independently interesting ideas. It is hard to build topological $\bfR^d$-bundles that do not arise from vector bundles, and essentially the only way to do so is as topological tangent bundles of topological manifolds, or more generally vertical topological tangent bundles of topological manifold bundles. The idea is to modify a smooth manifold bundle to produce a rather strange topological manifold bundle, using Disc-presheaves in the guise of the equivalence \eqref{eqn:bhomeo-as-embcalc} to give a pullback decomposition of $\BHomeo_\partial(M)$. Ideally, the latter would say that the sequence
	\[\BHomeo_\partial(M) \lra 
	\BAut_{\partial}(TM) \lra \BAut_{\partial}(M)\]
can be completed to a pullback square, but this statement needs modification. To see why, we explain how one obtains the pullback square: we start by using the smoothing theory of \cref{sec:smoothing-theory} to compare topological Disc-presheaves to ``particle'' Disc-presheaves as in \cref{sec:particle}, and the connection between them is controlled by the map $\BTop(d) \to \BAut(E_d)$. While we do not understand this map, we \emph{do} know that the map $\smash{\BTop(d-2) \to \BAut(E^\bfQ_d)}$ is null-homotopic, as used in \cref{sec:non-triviality}. Hence we add a suitable tangential structure $\theta$, fibrewise rationalisations $\fQ$, and restrictions to certain components---these impact the readability of the statement, but not its usefulness---to obtain the following result, in joint work of Krannich and myself \cite{KKSquare}:

\begin{thm}\label{thm:square} Suppose we have a compact topological $d$-dimensional triad $M$ for $d \geq 5$ with $M$ and $\partial_0 M$ both 2-connected and smoothable, and $\partial_1 M = D^{d-1}$. Let $\theta \colon \BSTop(d-2) \to \BSTop(d)$ and let $\ell_0$ be a $\theta$-structure on $\partial_0 M$. Then there is a space $Z_M$ and a pullback square
	\[\vspace{-.2cm}\begin{tikzcd} \BHomeo^\theta_\partial(M;\ell_0)_{\fQ} \rar \dar & Z_M \dar \\[-7pt]
		\BAut^\theta_{\partial}(M;\ell_0)_\fQ \rar & \BAut_{\partial}(M)^{\ell_0}.\end{tikzcd}\]
\end{thm}

How does one use \cref{thm:square} to manipulate a topological manifold bundle and prove \cref{thm:pontryagin-weiss}? Endowed with suitable tangential structure on its vertical tangent bundle, such a bundle gives a map into the top-left corner. By the universal property of a pullback this is the same data as (a) a fibration with vector bundle on its total space (the bottom-left corner), (b) whatever is classified by $Z_M$ (the top-right corner), (c) an isomorphism between their underlying fibrations (the homotopy in the bottom-right corner). This ``decoupling'' of the tangential data from the rest, allows one to set the generalised Miller--Morita--Mumford classes $\kappa_c$ for $c$ that are decomposable in the $p_i$ to zero, as long as one leaves $\kappa_{\scr{L}_i}$ fixed (by the family signature theorem these arise from the bottom-right corner). Starting with a suitable topological manifold bundle that has $\kappa_{\scr{L}_i} \neq 0$, this means there must also be a bundle that has $\kappa_{\scr{L}_i} \neq 0$ but all $\kappa_c = 0$ for $c$ that are decomposable in the $p_i$, so necessarily $\kappa_{p_i} \neq 0$ and hence $p_i \neq 0 \in H^{4i}(\BTop(d);\bfQ)$. Starting with a suitable smooth manifold bundle over a sphere, provided by the work of Galatius and Randal-Williams on stable moduli spaces, this argument yields in addition to the non-vanishing that the $p_i$ are detected on the Hurewicz image.

\subsection{Recovering boundaries and complements} Implicit in \cref{thm:square} is that one can, under certain conditions, recover the complement of the image of an embedding from the induced map of presheaves and hence the boundary as the complement of the interior.

\subsubsection{Recovering a complement} Suppose that $M \subset N$ is a compact codimension zero submanifold contained in the interior of $N$. Then the homotopy type of the complement  $N\backslash M$ can be obtained as follows: isotopy extension provides a fibre sequence
\[\Emb(\bfR^d,N \backslash M) \lra \Emb(\bfR^d \sqcup M,N) \lra \Emb(M,N)\]
with fibre taken over the inclusion, and evaluation at the origin induces an equivalence from the orbits $\smash{\Emb(\bfR^d,N \backslash M)_{\GL_d(\bfR)}}$ of the action by precomposition to $N \backslash M$. In analogy, we can define a \emph{$T_\infty$-complement} \cite{KKSquare}
\[N \backslash^{T_\infty} M \coloneqq \fib_\inc\big(T_\infty \Emb(\bfR^d \sqcup M,N) \to T_\infty \Emb(M,N) \big).\]
By the isotopy extension theorem for Disc-presheaves of \cref{thm:isotopy-extension} this is often equivalent to the actual complement, for instance if $M$ has handle dimension $\leq d-3$, an observation going back to Tillmann and Weiss \cite{TillmannWeiss}.

\subsubsection{Recovering the boundary} The above explains how it is possible that the right term $\rm{B}T_\infty \Emb^{t,\cong}_{\partial_0}(M,M)$ of \eqref{eqn:bhomeo-as-embcalc} seems only aware of the pair $(M,\partial_0 M)$, while the bottom-right entry of \cref{thm:square} involves the larger pair $(M,\partial M)$: Disc-presheaves recover the missing part of the boundary $\partial_1 M$, here $D^{d-1}$, as the complement of $\rm{int}(M) \cup \partial_0 M$ in $M$.

There is an analogue for moduli spaces. Suppose we have a topological manifold triad, given by a compact $d$-dimensional topological manifold $M$ whose boundary is divided into $\partial_0 M$ and $\partial_1 M$ meeting at $\partial_{01} M$. Any homeomorphism of $M$ that fixes $\partial_0 M$ pointwise must necessarily also preserve $\partial_1 M$ setwise, and hence there is a factorisation
\[\BHomeo_{\partial_0 M}(M) \lra \BAut_{\partial_0 M}(M,\partial_1 M) \lra \BAut_{\partial_0}(M)\]
of the map to homotopy automorphisms of $M$ fixing $\partial_0 M$ pointwise over those that in addition fix $\partial_1 M$ setwise. Extending work of Weiss \cite{WeissHaut}, Krannich and I proved the same is true for particle Disc-presheaves (as in \cref{sec:particle}, with attendant embedding calculus denoted $T_\infty \Emb^\rm{p}$) and thus for topological Disc-presheaves, under convergence hypotheses \cite{KKSquare}:

\begin{thm}If $d \geq 5$, $M$ is smoothable, and both inclusions $\partial_1 M \subset M$ and $\partial_{01} M \subset \partial_0 M$ induce equivalences on tangential 2-types, there is a dashed map making the following diagram commute
	\vspace{-.25cm}\[ \begin{tikzcd} \BHomeo_{\partial_1}(M) \rar \dar &[-20pt] \rar  \rm{B}T_\infty \Emb^\rm{t}_{\partial_0}(M,M)^\times &[-10pt] \rm{B}T_\infty \Emb^\rm{p}_{\partial_0}(M,M)^\times \dar  \arrow[dashed]{lld} \\[-7pt]
		\BAut_{\partial_0 M}(M,\partial_1 M) \arrow{rr} & & \BAut_{\partial_0}(M).\end{tikzcd}\]
\end{thm}

\vspace{-.25cm} \bibliographystyle{amsalpha}
\bibliography{literature}

\end{document}